\documentclass[journal]{IEEEtran}
\usepackage{xcolor,soul,framed} 

\colorlet{shadecolor}{yellow}
\usepackage[pdftex]{graphicx}

\usepackage[cmex10]{amsmath}
\usepackage{array}
\usepackage{url}
\usepackage{hyperref}

\usepackage{svg}
\usepackage{mathtools}

\usepackage{amssymb}

\newtheorem{assumption}{Assumption}
\newtheorem{theorem}{Theorem}
\newtheorem{definition}{Definition}
\newtheorem{remark}{Remark}
\newtheorem{example}{Example}

{\lccode`\?=`\p \lccode`\!=`\t  \lowercase{\gdef\ignorept#1?!{#1}}}
\def\usedim#1 {\expandafter\ignorept\the#1 \space}

\def\cyclicseq#1{\setbox0=\hbox{\kern-.7em$#1$\kern-.7em}%
    \dimen0=.3\wd0 \dimen1=.7\wd0
    \leavevmode \kern.7em
    \pdfsave\pdfsetmatrix{.9397 .342 -.342 .9397}\raise4pt\rlap{$\prec$}\pdfrestore
    \pdfliteral{q .9963 0 0 .9963 0 7 cm .4 w 0 0 m
        \usedim\dimen0 10 \usedim\dimen1 10 \usedim\wd0 0 c S Q}%
    \vbox to15pt{}\box0 \kern.7em
}

\begin{document}
    \title{Recursive Self-Composite Approach Towards Structural Understanding of Boolean Networks}
  \author{
      Jongrae~Kim,~\IEEEmembership{Senior Member,~IEEE,}
      Woojeong~Lee,~
      and
      Kwang-Hyun~Cho,~\IEEEmembership{Senior Member,~IEEE,}
      \\

  \thanks{Manuscript received  26 December 2023. 
  This work was supported by the Cheney Fellow program of the University of Leeds, U.K. It was also supported by the National Research Foundation of Korea (NRF) grants funded by the Korea Government, the Ministry of Science and ICT [2023R1A2C3002619 and 2021M3A9I4024447 (Bio \& Medical Technology Development Program)] and the internal fund/grant of Electronics and Telecommunications Research Institute (ETRI) [23RB1100, Exploratory and Strategic Research of ETRI-KAIST ICT Future Technology].}
  \thanks{J.~Kim is with the School of Mechanical Engineering, University of Leeds, Leeds LS2 9JT, UK (e-mail: menjkim@leeds.ac.uk).}
  \thanks{W.~Lee is with the Department of Bio \& Brain Engineering, KAIST, Daejeon, Republic of Korea (e-mail: frship35@kaist.ac.kr).}
  \thanks{K-H.~Cho is with the Department of Bio \& Brain Engineering, KAIST, Daejeon, Republic of Korea (e-mail: ckh@kaist.ac.kr, Corresponding author).}
  }  


\twocolumn[{%
This work has been submitted to the IEEE for possible publication. Copyright may be transferred without notice, after which this version may no longer be accessible.
}]

\maketitle

\begin{abstract}
Boolean networks have been widely used in many areas of science and engineering to represent various dynamical behaviour. In systems biology, they became useful tools to study the dynamical characteristics of large-scale biomolecular networks and there have been a number of studies to develop efficient ways of finding steady states or cycles of Boolean network models. On the other hand, there has been little attention to analyzing the dynamic properties of the network structure itself. Here, we present a systematic way to study such properties by introducing a recursive self-composite of the logic update rules. Of note, we found that all Boolean update rules actually have repeated logic structures underneath. This repeated nature of Boolean networks reveals interesting algebraic properties embedded in the networks. We found that each converged logic leads to the same states, called kernel states. As a result, the longest-length period of states cycle turns out to be equal to the number of converged logics in the logic cycle. Based on this, we propose a leaping and filling algorithm to avoid any possible large string explosions during the self-composition procedures. Finally, we demonstrate how the proposed approach can be used to reveal interesting hidden properties using Boolean network examples of a simple network with a long feedback structure, a T-cell receptor network and a cancer network.
\end{abstract}

\begin{IEEEkeywords}
Boolean networks, logic structures, kernel states, biological networks, systems biology
\end{IEEEkeywords}

\section{Introduction}

\IEEEPARstart{B}{oolean} network formalism is a useful mathematical modelling approach to describe complex interactions and dynamics of biological systems \cite{kauffman1969homeostasis}. In this formalism, individual biological entities, such as genes, proteins, or other molecular components, are represented by nodes, while their interactions are depicted as edges. These nodes are assigned with time-varying binary states – either on (active) or off (inactive) – thus facilitating a simplified modelling process and allowing for a broader range of interactions while still capturing essential dynamical properties. Logical relationships among these nodes are specified through Boolean functions.
Following these rules, node states are updated synchronously or asynchronously, eventually converging to a stable state known as an attractor. It has been previously demonstrated that stable attractor states in gene regulatory networks correspond to distinct cellular phenotypes or cell fates \cite{huang2005cell}.
In this regard, extensive studies have been done to investigate the long-term behaviour of biological networks represented by Boolean network models, aiming to predict real intra-cellular dynamics across various biological processes, including the cell cycle \cite{sahin2009modeling}, differentiation \cite{collombet2017logical, carbo2013systems}, and tumorigenesis \cite{fumia2013boolean, cho2016attractor}.
Such studies have not only enhanced our understanding of biological phenomena but also enabled the prediction of drug responses for precision medicine of complex diseases \cite{choi2017network} and the identification of potential therapeutic targets for drug discovery \cite{choi2022network, kim2023cell}.
\bigskip

Let us consider the Boolean networks given by
\begin{align} \label{eq:boolean_network_1_to_n}
    x_1(k+1) &= f_1[x_1(k), x_2(k), \ldots, x_{n-1}(k), x_n(k)]\notag\\
    x_2(k+1) &= f_2[x_1(k), x_2(k), \ldots, x_{n-1}(k), x_n(k)]\notag\\
    &\vdots\notag\\
    x_n(k+1) &= f_n[x_1(k), x_2(k), \ldots, x_{n-1}(k), x_n(k)]
\end{align}
where 
$x_i(k)$ is the $i$-th Boolean state equal to either \emph{true} (equivalently T or 1) or \emph{false} (equivalently F or 0) at $k$ for $i=1, 2, \ldots, n-1, n$, 
$k$ is the non-negative integer in $[0, \infty)$, 
$x_i(0)$ is the initial state, 
$f_i(\cdot)$ is a synchronous update rule consisting of 
the Boolean operations conjunction (and, $\wedge$), disjunction (or, $\vee$), 
and negation(not, $\neg$)
and
$x_i(k+1)$ is the updated Boolean state for $i=1, 2, \ldots, n-1, n$.

The Boolean network shown in \eqref{eq:boolean_network_1_to_n} can be written in a compact form as follows:
\begin{align} \label{eq:boolean_network_compact}
    {\bf x}(k+1) = {\bf f}[{\bf x}(k)]
\end{align}
where
\begin{subequations}
\begin{align}
    {\bf x}(k) &= \begin{bmatrix} 
                x_1(k) & x_2(k) & \ldots & x_n(k)
                \end{bmatrix}^T\\
    {\bf f}[{\bf x}(k)] &= \begin{bmatrix} 
                f_1[{\bf x}(k)] & f_2[{\bf x}(k)] & \ldots & 
                f_n[{\bf x}(k)]^T
                \end{bmatrix}
\end{align}
\end{subequations}
and $(\cdot)^T$ is the transpose.

The state space is $2^n$-dimensional and the main interest in Boolean network analysis is finding steady-states and periodic cycles. In the synchronous update of Boolean networks, every initial state converges to a steady state or a periodic cycle. As $n$ increases, the dimension of the state space, $2^n$, increases exponentially. Therefore, executing the exhaustive search to find attractors is infeasible even for moderate-size networks, e.g., $n$ around 30. One of the well-known approaches in Boolean networks called the semi-tensor approach is also an exhaustive method \cite{cheng2010linear}. The aggregation algorithm proposed in \cite{6475982} relies on the specific modular structure of the networks. Hence, Boolean network analysis results are often obtained from probabilistic approaches based on simulations over a finite number of random samples. Finding attractors or control strategies for Boolean networks is also known to be NP-hard \cite{akutsu2007control}.

\begin{figure}
    \centering
    \includegraphics[width=0.45\textwidth]{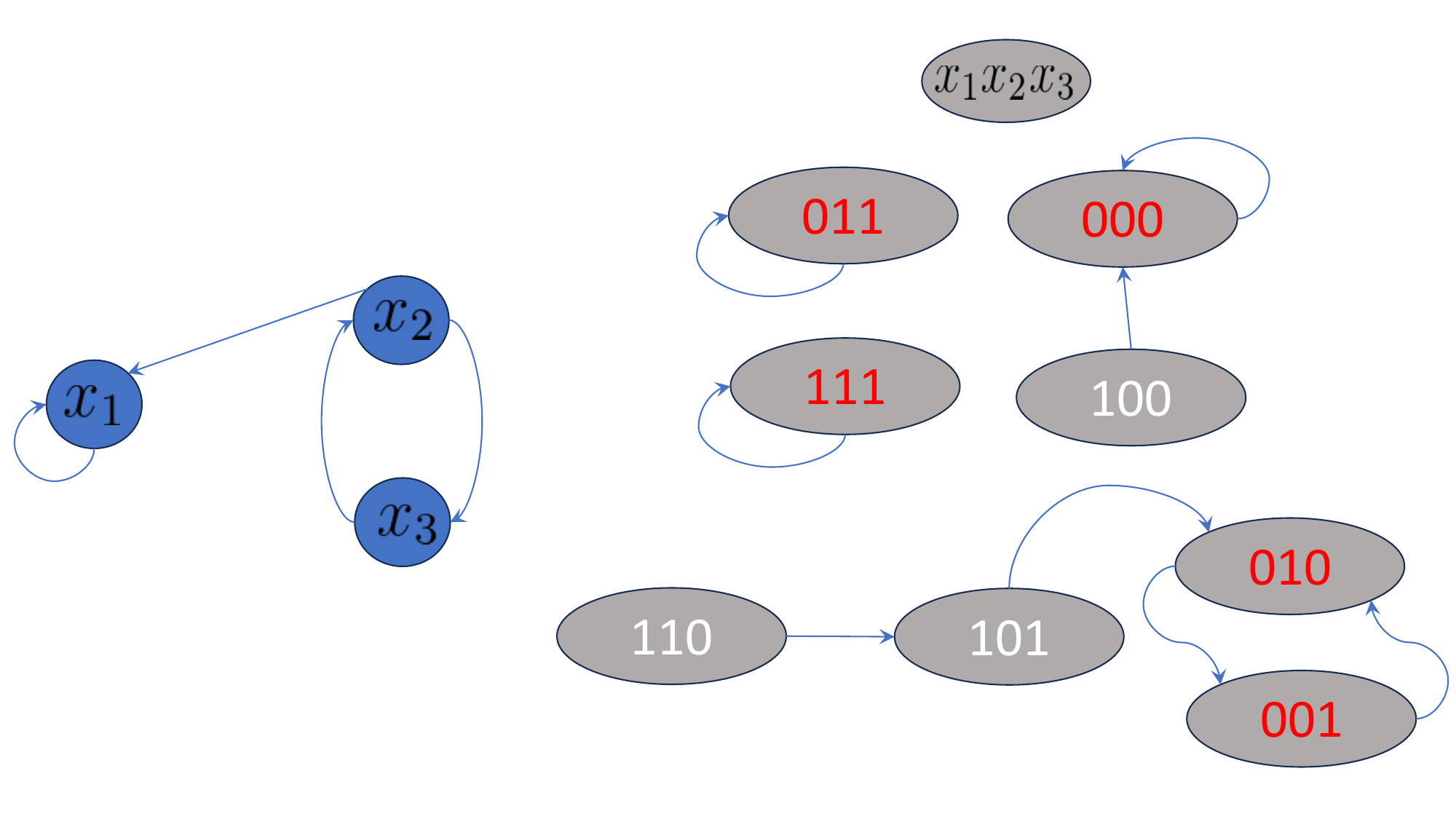}
    \caption{The interaction graph and state transition map of the Boolean network, \eqref{eq:x123_simple_example_01}, are shown.  There are three steady-states, $\{ s = (000), (011), (111) \}$ and one cycle with the period 2, $\{ s= (010) \leftrightarrow (001) \}$, where $s=(x_1 x_2 x_3)$.}
    \label{fig:x123_BL}
\end{figure}

On the other hand, there are rich theoretical results in the continuous system described by ordinary differential equations: $dx/dt = f(x)$, where $d(\cdot)/dt$ is the derivative with respect to time, $t$, and $f(x)$ is a nonlinear function satisfying the existence and uniqueness conditions of the solution. The equilibrium points and their stability are inherent properties of the right-hand side of the differential equation, i.e., $f(x)$. The solution of $f(x) = 0$ is the equilibrium point and the eigenvalues of $df(x)/dx$ at the equilibrium point provide the stability condition. 
The main motivation of our approach to be shown is from the question about whether the right-hand side of \eqref{eq:boolean_network_compact}, i.e., ${\bf f}[{\bf x}(k)]$, can also provide any clue for the characteristics of Boolean networks.

In the following sections, first, the motivation of the proposed method is illustrated with a simple toy example. Secondly, we present the main results of the recursive self-composition approach to investigate the structure of Boolean networks. Thirdly, we apply the proposed method to various examples including a simple network with a long feedback path and two biological networks – a T-cell signalling pathway and a cancer signalling network – highlighting the advantages of the proposed method. Finally, the conclusions are made.

\section{Recursive Self-Composite Boolean network}

\subsection{Motivations}
Let us consider the following Boolean network model:
\begin{subequations} \label{eq:x123_simple_example_01}
\begin{align}
    x_1(k+1) &= x_1(k) \wedge x_2(k)\\
    x_2(k+1) &= x_3(k)\\
    x_3(k+1) &= x_2(k)
\end{align}
\end{subequations}
Since the right-hand side of $x_1(k+1)$ in \eqref{eq:x123_simple_example_01}  is of the conjunction of $x_1(k)$ and $x_2(k)$, the 75\% of $x_1(k+1)$ is 0 (False), i.e., the 25\% of $x_1(k+1)$ is 1 (True). All four possible outputs from the conjunction of $x_1(k)$ and $x_2(k)$ produce 0 except when both are 1. Using the same approach, examining the right-hand sides of $x_2(k+1)$ and $x_3(k+1)$ in \eqref{eq:x123_simple_example_01}, the probability that the output of $x_2(k+1)$ or $x_3(k+1)$ is 0 or 1 is 50\%.

The question is how accurate these probabilities are with respect to the final state. The final state is a steady state or a state belonging to a cycle. Figure \ref{fig:x123_BL} shows the interaction graph and transition map, where the state, $s$, is equal to $(x_1 x_2 x_3)$. Three steady states and one cycle are shown in red. If we consider all eight states in Figure \ref{fig:x123_BL} and count the number of states with the final state $x_1=1$, there is only one case. In all other cases except $(x_1 x_2 x_3) = (1 1 1)$, $x_1$ becomes 0. Hence, the probability for $x_1$ equal to 0 in the final state is 87.5\%, 7 out of 8. It is not equal to the 75\% that was estimated earlier.

The cause for this difference is the usage of one-step propagation equation. So, a longer propagation would provide better estimation. Any exact calculation can be done by exhaustive numerical simulation considering all possible states, which is not feasible for large-size networks. Instead, let us consider the two-step propagation symbolically as follows:
\begin{subequations}
\begin{align}
    x_1(k+2) &= x_1(k+1) \wedge x_2(k+1) \notag\\
             &= \left[ x_1(k) \wedge x_2 (k) \right] \wedge x_3(k) \notag\\
             &= x_1(k) \wedge x_2 (k) \wedge x_3(k)
              \label{eq:x123_simple_example_01_two_step}\\
    x_2(k+2) &= x_3(k+1) = x_2(k)\\
    x_3(k+2) &= x_2(k+1) = x_3(k)
\end{align}
\end{subequations}
Similarly, the three-step propagation is obtained as follows:
\begin{subequations}    
\begin{align}
    x_1(k+3) &= x_1(k+1) \wedge x_2 (k+1) \wedge x_3(k+1)\notag\\
             & = \left[ x_1(k) \wedge x_2(k) \right]
             \wedge x_3(k) \wedge x_2(k)\notag\\
             &= x_1(k) \wedge x_2 (k) \wedge x_3(k)
              \label{eq:x123_simple_example_01_three_step}\\
    x_2(k+3) &= x_2(k+1) = x_3(k)\\
    x_3(k+3) &= x_3(k+1) = x_2(k)
\end{align}
\end{subequations}
The four-step propagation is given by
\begin{subequations}    
\begin{align}
    x_1(k+4) &= x_1(k+1) \wedge x_2 (k+1) \wedge x_3(k+1)\notag\\
             &= x_1(k) \wedge x_2 (k) \wedge x_3(k)
              \label{eq:x123_simple_example_01_four_step}\\
    x_2(k+4) &= x_3(k+1) = x_2(k)\\
    x_3(k+4) &= x_2(k+1) = x_3(k)
\end{align}
\end{subequations}
and it turns out that the five-step propagation is the same as the two-step propagation. 

\begin{figure}
    \centering
    \includegraphics[width=0.47\textwidth]{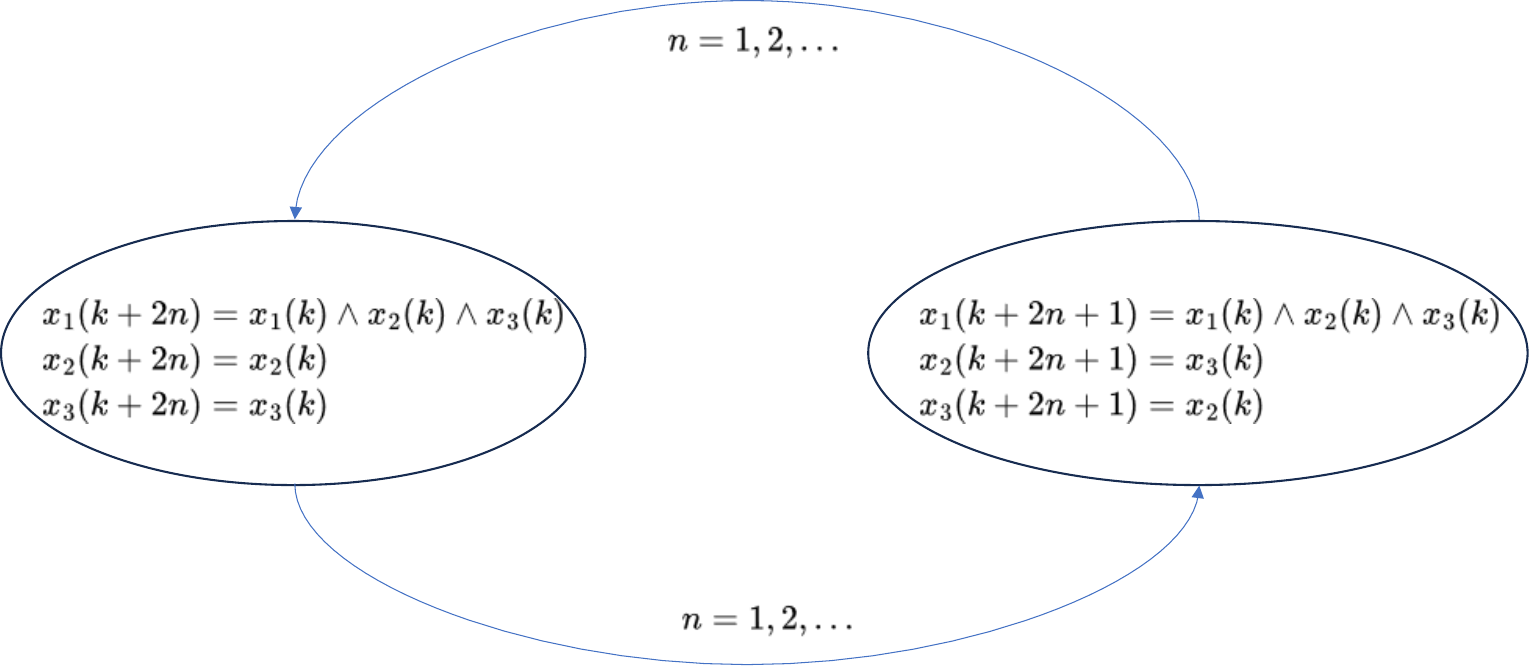}
    \caption{The Boolean logic in \eqref{eq:x123_simple_example_01} switches between two update rules.}
    \label{fig:x1x2x3_logic_cycle}
\end{figure}

As shown in Figure \ref{fig:x1x2x3_logic_cycle}, the update logic itself switches between the two update rules. While the $x_2$ and $x_3$ propagation rules cross-update between the two, the $x_1$ update rule converges to the conjunction of the three states.
By inspecting the right-hand side of the converged update rule, the probability of $x_1$ converging to the final state equal to 0 is 7 out of 8, which coincides with the true probability.

In the exhaustive approach, the transition matrix, $L$, describes the updates of eight states in Figure \ref{fig:x123_BL} as follows:
\begin{align}
    {\bf s}_{k+1}
    = 
    \underbrace{
    \begin{bmatrix}
    1 & 0 & 0 & 0 & 0 & 0 & 0 & 0\\
    0 & 0 & 1 & 0 & 0 & 0 & 0 & 0\\
    0 & 1 & 0 & 0 & 0 & 0 & 0 & 0\\
    0 & 0 & 0 & 1 & 0 & 0 & 0 & 0\\
    1 & 0 & 0 & 0 & 0 & 0 & 0 & 0\\
    0 & 0 & 1 & 0 & 0 & 0 & 0 & 0\\
    0 & 0 & 0 & 0 & 0 & 1 & 0 & 0\\
    0 & 0 & 0 & 0 & 0 & 0 & 0 & 1\\
    \end{bmatrix}
    }_{L}
    \underbrace{
    \begin{bmatrix}
    (000)\\ (001)\\ (010)\\ (011)\\ (100)\\ (101)\\ (110)\\ (111)
    \end{bmatrix}
    }_{{\bf s}_{k}}
\end{align}
In the semi-tensor approach \cite{cheng2010linear}, ${\bf s}_k$ or ${\bf s}_{k+1}$ is the state vector
with the element corresponding to the current state equal to 1 and the rest set to 0, i.e.,
\begin{align}
    {\bf s}_k = \begin{bmatrix}
        s_0 & s_1 & \ldots & s_{2^n-1}
    \end{bmatrix}^T
\end{align}
where $s_i = 1$ for $i$ equal to the decimal number whose binary number
corresponds to the current state $(x_1 x_2 \ldots x_n)$ and $s_i = 0$ for the others,
where $i \in \{0, 1, 2, \ldots, 2^n-1\}$.
For instance, if the initial state for $(x_1, x_2, x_3)$ is equal to (010), then, $s_2 = 1$ and
the rest of $s_i$ equal to 0.
\begin{align}
    {\bf s}_0 = \begin{bmatrix}
        0 & 0 & 1 & 0 & 0 & 0 & 0 & 0
    \end{bmatrix}^T
\end{align}
(010), i.e., $s_2=1$, is converted to (001), i.e., $s_1=1$, as shown in Figure \ref{fig:x123_BL}, and $L {\bf s}_0$ provides the corresponding transition state, ${\bf s}_1$, i.e.,
\begin{align}
    {\bf s}_{1} = \begin{bmatrix}
        0 & 1 & 0 & 0 & 0 & 0 & 0 & 0
    \end{bmatrix}^T
\end{align}
Further iterations reveal that $L^{2k} = L^2$ and $L^{2k+1} = L^3$ for $k=1, 2, \ldots$. Hence, the algorithms constructed in \cite{cheng2010linear} can find all steady states and cycles by inspecting $L^2$ and $L^3$. However, one of the main drawbacks of this approach is the requirement for \emph{always} checking all $2^n$ states to construct $L$. Hence, the algorithm is limited to solving only small or moderate-size Boolean networks only.

On the other hand, the approach we propose does not require explicitly checking $2^n$ states or constructing the matrix $L$. As shown in Figure \ref{fig:x1x2x3_logic_cycle}, the recursive self-composition provides the converged cyclic logic without checking the $2^3$ states.

\subsection{Main Results}

\begin{assumption}[Synchronous Update]
All states in the Boolean network given by \eqref{eq:boolean_network_compact} are updated synchronously.
The states in the right-hand-side of \eqref{eq:boolean_network_compact} is the one at the same step.
\end{assumption}
\bigskip

\begin{definition}[Recursive Self-Composite] \label{theorem:recur_self_composite}
    The $p$-times recursive self-composite of the Boolean network is given by
    \begin{align}
        {\bf x}(p) &= {\bf f}[{\bf x}(p-1)] = {\bf f}[{\bf f}({\bf x}(p-2))] = \ldots \notag\\
         &= {\bf f}[{\bf f}({\bf f}\ldots ({\bf f}({\bf x}(0)))] \notag\\
         &= \underbrace{{\bf f}\circ{\bf f}\circ{\bf f}\circ\ldots\circ{\bf f}}_{p-\text{times}}[{\bf x}(0)] 
         = {\bf f}^p[{\bf x}(0)]
    \end{align}
    where $p$ is a positive integer.
\end{definition}
\bigskip

\begin{theorem}[Convergence of Recursive Self-Composition]
    All recursive self-composition of the synchronous Boolean network given by
    \eqref{eq:boolean_network_compact} converges to a steady-state logic or a cyclic logic, i.e,
    \begin{align} \label{eq:converged_logic}
        {\bf f}^{p^* + r}({\bf x}) = {\bf f}^{p^* + r + \ell^*}({\bf x})
    \end{align}
where $p^*$ between 0 and $2^n$ is the minimum number of recursions when the logic starts repeating itself,
$r$ is a non-negative integer and
$\ell^*$ between 1 and $2^n$ is the period of logic cycles.
\end{theorem}
\emph{Proof:} Deterministic synchronous update Boolean networks have a finite number of states, i.e., $2^n$,
any initial state repeats the same state at longest in $2^n+1$ steps. 
It is guaranteed at least ${\bf x}(2n+1)$ is equal to ${\bf x}(p)$ for $1 \le p \le 2^n$.
As ${\bf x}(2n+1) = {\bf f}^{2n+1}[{\bf x}(0)]$ and ${\bf x}(p) = {\bf f}^p[{\bf x}(0)]$,
${\bf f}^p[{\bf x}(0)] = {\bf f}^{2n+1}[{\bf x}(0)]$ and $p$ is between 1 and $2^n$.
Hence, $p^*$ always exists between 1 and $2^n$. 

Assuming that the smallest $\ell^*$ is strictly greater than $2^n$ leads to a contradiction with the finite number of states $2^n$. Therefore, $\ell^*$ must be between 0 (for the cases of no periodic cycles but only steady states) 
and $2^n$. $\blacksquare$
\bigskip

When $\ell^*$ is equal to 0, it corresponds to the case that there is only one steady-state logic. 
If $\ell^*$ is equal to 3, three update rules switch between them.
For instance, the Boolean network given in \eqref{eq:x123_simple_example_01} has $p^*$ equal to 2, 
when the logic repetition begins to start, and $\ell^*$ equal to 2, which is the period of logic cycles. 

\bigskip

\begin{remark}
    In the worst case, if $p^*$ is equal to $2^n$, the computational cost to find a repeating logic
    is at least as much expensive as the exhaustive search.
\end{remark}

\bigskip

\begin{theorem}[Longest Cycle Upper Bound]
All cycle lengths of every Boolean logic cannot be longer than 
the length of the converged logic cycle, $\ell^*$.
\end{theorem}
\emph{Proof:}
Let us assume there exists a cycle of the period, $m$, strictly longer than $\ell^*$.
Let ${\bf s}_{p^*}$ be the state for the first time arrived in the cycle at $p^*$-step 
from the state ${\bf s}_{p^*-1}$, which is not in the cycle, and the cycle propagates as follows:
\begin{align*}
\begin{array}{ccccccc}
    {\bf s}_{p^*-1} & \xrightarrow{L} &  {\bf s}_{p^*} & \xrightarrow{L} & {\bf s}_{p^*+1} &  
    \xrightarrow{L} & {\bf s}_{p^*+2}\\
    & & & & & & \downarrow\text{\scriptsize $L$}\\
    {\bf s}_{p^*+\ell^*} & \xleftarrow{L} 
    & {\bf s}_{p^*+\ell^*+1} & \xleftarrow{L} & \ldots & \xleftarrow{L} & {\bf s}_{p^*+3}\\
    \text{\scriptsize $L$}\downarrow & & & & & \\
    {\bf s}_{p^*+\ell^*+1} & \xrightarrow{L} &  \ldots & \xrightarrow{L} & {\bf s}_{p^*+m-1} &  
    \xrightarrow{L} & {\bf s}_{p^*}
\end{array}
\end{align*}
where ${\bf s}_{p^*}$ in the last line is equal to the one in the first line and the cycle repeats.
By the definition of cycle, ${\bf s}_i \neq {\bf s}_j$ for $i \neq j$. 
Also, notice that as 
${\bf s}_{p^*+m} = {\bf s}_{p^*}$ and ${\bf f}^{p^*+m}({\bf x}_0) = {\bf f}^{p^*}({\bf x}_0) = {\bf x}_0$, 
$m$ must be an integer multiple of $\ell^*$.

We can re-write $(p^*+\ell^*+1)$-step using the composition logic as follows:
\begin{align*}
 &{\bf s}_{p^*+\ell^*+1} = L {\bf s}_{p^*+\ell^*}
 \rightarrow {\bf x}(p^*+\ell^*+1) = {\bf f}^{p^*+\ell^*+1}[{\bf x}(0)]
\end{align*}
where ${\bf x}(p^*+\ell^*+1)$ corresponds to ${\bf s}_{p^*+\ell^*+1}$.
By Theorem \ref{theorem:recur_self_composite}, the following equality satisfies
\begin{align}
{\bf f}^{p^*+\ell^*+1}[{\bf x}(0)] = {\bf f}^{p^*+1}[{\bf x}(0)]    
\end{align}
Hence,
\begin{align}
{\bf s}_{p^*+\ell^*+1} = {\bf s}_{p^*+1}    
\end{align}
This contradicts ${\bf s}_i \neq {\bf s}_j$ for $i \neq j$ in the cycle. Therefore,
no cycle can have a longer period than $\ell^*$. $\blacksquare$
\bigskip

\begin{definition}[Kernel States Set]
    The kernel states set, ${\mathbb K}$, of the synchronous Boolean network,
    \eqref{eq:boolean_network_compact},
    includes all steady-states and the states belonging to cycles.
\end{definition}
\bigskip

For instance, all the states indicated in red in Figure \ref{fig:x123_BL}
are the kernel states of the network and 
${\mathbb K} = \{(000), (001), (010), (011), (111) \}$
or equivalently ${\mathbb K} = \{ s_0, s_1, s_2, s_3, s_7 \}$.
\bigskip

\begin{theorem}[Kernel States Set of Converged Logic] \label{theorem:K_set_converged_logic}
    The one-step propagated state, ${\bf x}(1)$, by the converged logic 
    in \eqref{eq:converged_logic}, i.e., 
    ${\bf x}(1) = {\bf f}^{p^*+r}[{\bf x}(0)]$,
    converges
    the same kernel states set of the original Boolean network
    for any fixed non-negative integer $r$.
\end{theorem}
\emph{Proof:} If a steady-state is absent in the kernel set of a converged logic,
then it contradicts the property of steady-states. Hence, the proof for
steady-states cases becomes trivial.
Let us consider a cycle of the period, $\ell^*$, as follows:
\begin{align}
    \cyclicseq{{\bf s}_0 \underbrace{\to}_{L {\bf s}_0} {\bf s}_1 \underbrace{\to}_{L {\bf s}_1} 
    {\bf s}_2 \underbrace{\to}_{L{\bf s}_2} \ldots \underbrace{\to}_{L{\bf s}_{\ell^*-1}} {\bf s}_{\ell^*}} 
\end{align}
Let the initial state, ${\bf x}(0)$, correspond to ${\bf s}_0$. And, it propagates
$p^*$ steps as follows:
\begin{align}
    {\bf x}_1 = {\bf f}[{\bf x}(0)] &\to {\bf x}_2 = {\bf f}^1[{\bf x}(1)]
    \to \ldots \notag\\
    \ldots &\to {\bf x}_{p^*} = {\bf f}^{p^*}[{\bf x}(p^*-1)] 
\end{align}
As ${\bf x}(0)$ starts in the cycle, all propagated states are in the cycle. 
Hence, ${\bf x}_{p^*}$ is equal to one of the states in the cycle. Specifically,
${\bf x}_{p^*}$ is equal to the state corresponding to ${\bf s}_{r^*}$,
where $r^* = p^* - \ell^* q$, which is in $[0, \ell^*]$, 
$q$ is the maximum integer such that $\ell^* q$ is less
than or equal to $p^*$. 

Without loss of the generality, let us assume that $r^*$ is equal to 2, i.e.,
${\bf x}_p^*$ corresponds to ${\bf s}_2$. It implies that:
${\bf s}_2$ is an element of the kernel state set of ${\bf f}^{p^*}$,
${\bf s}_3$ is an element of the kernel state set of ${\bf f}^{p^*+1}$
and so forth.

Let us choose the initial state, ${\bf x}(0)$ corresponding to ${\bf s}_1$ and
repeat the same procedure. Then, $r^*$ becomes 3 and this results in: 
${\bf s}_3$ is an element of the kernel state set of ${\bf f}^{p^*}$,
${\bf s}_4$ is an element of the kernel state set of ${\bf f}^{p^*+1}$
and so forth.

For the shorter cycles less than the period $\ell^*$, the same steps
provide the proof that all the cycle states must be in the kernel state
of each of the converged logic.
Therefore, the converged logic includes all states in the cycles. $\blacksquare$
\bigskip

\begin{theorem}[Longest Length Cycle]
    There exists at least one cycle whose length is equal to the period
    of the logic cycle, $\ell^*$.
\end{theorem}
\emph{Proof:} By Theorem \ref{theorem:K_set_converged_logic}, the range
set of every converged logic is identical with each other as the kernel set, ${\mathbb K}$. 
Hence, once the logic converges to the logic cycle, whose period is $\ell^*$,
the mapping from ${\mathbb K}$ to ${\mathbb K}$ repeats $\ell^*$ times.
Each of the mappings must be different from each other. Otherwise,
the existence of the logic cycle equal to $\ell^*$ is violated.
In addition, due to the periodicity of the logic cycles, 
the $\ell^*$-th mapping brings the states back to the states mapped by
the first logic cycle. $\blacksquare$
\bigskip

\begin{definition}[Kernel Logic]
    The Kernel logic, ${\bf f}^{p^*+k^*}$ is the converged logic having the same steady states and cycles as the original
    Boolean logic, where $k^*$ is an integer between 0 and $\ell^*-1$.
\end{definition}
\bigskip
For instance, the converged logic in the right-hand side of Figure \ref{fig:x1x2x3_logic_cycle} has
the same three steady states and one cycle as the original network given by \eqref{eq:x123_simple_example_01}.
\bigskip

\begin{figure*}
    \centering
    \includegraphics[width=0.87\textwidth]{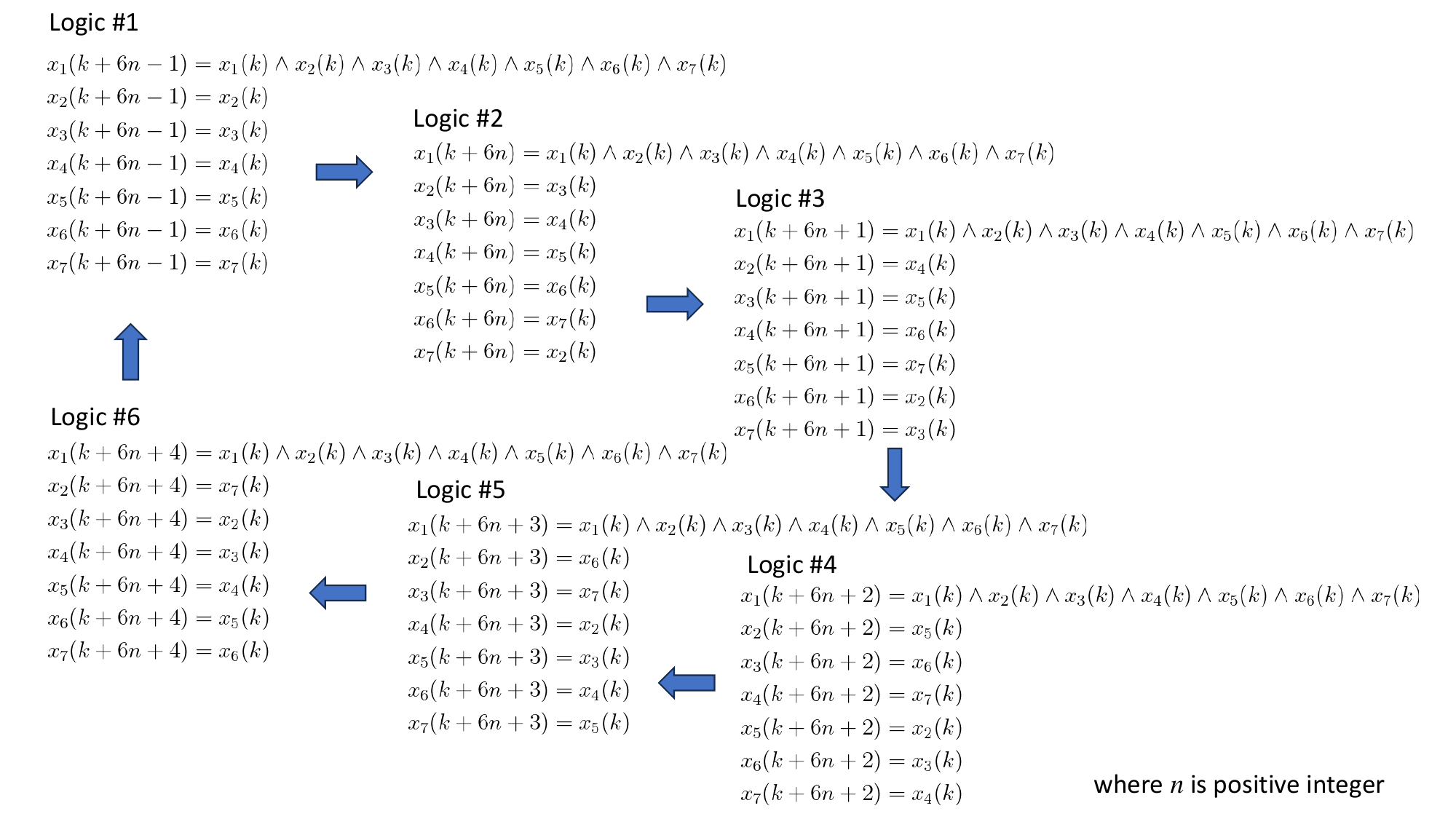}
    \caption{The Boolean logic in \eqref{eq:x1234567_simple_example} cycles between six update rules.}
    \label{fig:x1x2x3x4x5x6x7_logic_cycle}
\end{figure*}

\begin{theorem}[Existence of Kernel Logic]
For every Boolean network, there exists at least one kernel logic among the converged logic cycles.
\end{theorem}
\emph{Proof:} Given that there are $\ell^*$ number of converged logic, ${\bf f}^{p^*+k}$, where
$k$ is an integer from 0 to $\ell^*-1$, take $n$-time self composites for a fixed $k$ logic
as follows:
\begin{align}
    \underbrace{{\bf f}^{p^*+k}\circ \ldots \circ {\bf f}^{p^*+k} \circ {\bf f}^{p^*+k}}_{t-\text{times}}[{\bf x}(0)]
    = {\bf f}^{t(p^*+k)}[{\bf x}(0)]
\end{align}
$p^*+k$ can be expressed as
\begin{align}
    p^*+k = m \ell^* + r
\end{align}
where $m$ is the largest integer satisfying $p^*+k \ge m \ell^*$, where $r$ can be any integer between 0 and 
$\ell^*-1$ as $k$ is between 0 and $\ell^*-1$. 
Set $r=1$ and multiply $t$, which is an integer,
\begin{align}
    t(p^*+k) = t m \ell^* + t
\end{align}
Hence, we can cover all integers from 0 to $\ell^*-1$
by varing $t$. Therefore, it covers all $\ell^*$ cyclic
logic and there exists always at least one kernel logic. $\blacksquare$
\bigskip

\begin{remark}[Leaping \& Filling]
One of the ways to speed up the self-composition iteration and possibly increase
the chance to avoid large string-length explosions is leaping by performing larger-step
composition instead of a one-step composition.
First, 
\begin{align*}
{\bf x}(k+2) = {\bf f}^2[{\bf x}(k)] = {\bf g}[{\bf x}(k)]
\end{align*}
is obtained. Secondly, 
\begin{align*}
    {\bf x}(k+4) = {\bf g}[{\bf x}(k+2)] = {\bf g}^2[{\bf x}(k)] = {\bf h}[{\bf x}(k)]
\end{align*}
then, 
\begin{align*}
    {\bf x}(k+8) = {\bf h}[{\bf x}(k+4)] = {\bf h}^2[{\bf x}(k)]
\end{align*}
and we continue until the logic converges. Once the logic converges, we apply ${\bf f}(\cdot)$
repeatedly and obtain the logic between the leaps. For instance, ${\bf h}[{\bf x}(k)]$ converges
and ${\bf h}^2[{\bf x}(k)]$ is equal to ${\bf h}[{\bf x}(k)]$. Then, the logic for
${\bf x}(k+5)$, ${\bf x}(k+6)$ and ${\bf x}(k+7)$ are obtained by the filling sequence as follows:
\begin{align*}
    {\bf x}(k+5) &= {\bf f}[{\bf x}(k+4)] = {\bf f}\{{\bf h}[{\bf x}(k)]\}\\
    {\bf x}(k+6) &= {\bf f}[{\bf x}(k+5)] = {\bf f}^2\{{\bf h}[{\bf x}(k)]\}\\
    {\bf x}(k+7) &= {\bf f}[{\bf x}(k+6)] = {\bf f}^3\{{\bf h}[{\bf x}(k)]\}
\end{align*}
\end{remark}

\section{Examples}

\begin{example}[Longer feedback network]
    The following Boolean network is an extended network of \eqref{eq:x123_simple_example_01}
    having a longer feedback chain between $x_2$ and $x_7$:
    \begin{subequations} \label{eq:x1234567_simple_example}        
    \begin{align} 
        x_1(k+1) &= x_1(k) \wedge x_2(k)\\
        x_2(k+1) &= x_3(k)\\
        x_3(k+1) &= x_4(k)\\
        x_4(k+1) &= x_5(k)\\
        x_5(k+1) &= x_6(k)\\
        x_6(k+1) &= x_7(k)\\
        x_7(k+1) &= x_2(k)
    \end{align}
    \end{subequations}
    By the recursive compositions, the logic converges at $p^* = 5$ to the cyclic logic whose period, 
    $\ell^*$, is equal to 6 as shown in Figure \ref{fig:x1x2x3x4x5x6x7_logic_cycle}. 
    It has three steady-states, a cycle with period 2, two cycles with period 3 
    and nine cycles with period 6. 
    The number of elements in the kernel states set is 65 ($=3 + 1\times 2 + 2\times 3 + 9\times 6$). 
    Among the six logics in the cycle, Logic \#2 and \#6 in Figure \ref{fig:x1x2x3x4x5x6x7_logic_cycle} 
    are the kernel logic.
    By inspecting the right-hand side of the converged logic, 
    the common characteristics found are as follows: 
    $x_1$ equal to 1 is a rare event, and the other states equal to 0 or 1 have the same chance.
    This coincides with the fact that there is only one kernel state with $x_1$ equal to 1 among
    the 65 kernel states.
    Hence, we would be able to design numerical
    and lab experiments to find rare events, which would be challenging to find even in numerical
    experiments based on Monte Carlo type random simulations.
\end{example}

\bigskip
\begin{example}[T-cell receptor network]
    A T-cell receptor Boolean network model in \cite{klamt2006methodology}
    having 37 states with 3 control inputs is given by
    \begin{align*}
    \text{AP1}(k+1) &= \text{Fos}(k) \wedge \text{Jun}(k),~
    \text{Ca}(k+1) = \text{IP3}(k)\\  
    \text{Calcin}(k+1) &= \text{Ca}(k),~
    \text{cCbl}(k+1) = \text{ZAP70}(k)\\
    \text{CRE}(k+1) &= \text{CREB}(k),~
    \text{CREB}(k+1) = \text{Rsk}(k)\\         
    \text{DAG}(k+1) &= \text{PLCg$^*$}(k),~
    \text{ERK}(k+1) = \text{MEK}(k)\\
    \text{Fos}(k+1) &= \text{ERK}(k)\\
    \text{Fyn}(k+1) &= [\text{Lck}(k) \wedge \text{CD45}(k)]\\
             &\vee [\text{TCR$^+$}(k) \wedge \text{CD45}(k)]\\
    \text{Gads}(k+1) &= \text{LAT}(k),~ 
    \text{Grb2Sos}(k+1) = \text{LAT}(k)\\
    \text{IKKbeta}(k+1) &= \text{PKCth}(k),~
    \text{IP3}(k+1) = \text{PLCg(act)}(k)\\
    \text{Itk}(k+1) &= \text{SLP76}(k) \wedge \text{ZAP70}(k)\\
    \text{IkB}(k+1) &= \neg \text{IKKbeta}(k),~
    \text{JNK}(k+1) = \text{SEK}(k)\\
    \text{Jun}(k+1) &= \text{JNK}(k),~
    \text{LAT}(k+1) = \text{ZAP70}(k)\\
    \text{Lck}(k+1) &= \neg \text{PAGCsk}(k) \wedge \text{CD45}(k) \wedge \text{CD4}(k)\\
    \text{MEK}(k+1) &= \text{Raf}(k),~
    \text{NFAT}(k+1) = \text{Calcin}(k)\\
    \text{NFkB}(k+1) &= \neg \text{IkB}(k),~
    \text{PKCth}(k+1) = \text{DAG}(k)\\
    \text{PLCg$^*$}(k+1) &=  [\text{Itk}(k) \wedge \text{PLCg$^+$}(k) \wedge \text{SLP76}(k)\\
                    &\wedge \text{ZAP70}(k)]
                    \vee [\text{PLCg$^+$}(k) \wedge \text{Rlk}(k)\\ 
                    &\wedge \text{SLP76}(k) \wedge \text{ZAP70}(k)]
    \end{align*}
    \begin{align*}                    
    \text{PAGCsk}(k+1) &= \text{Fyn}(k) \vee \neg \text{TCR$^+$}(k)\\
    \text{PLGg$^+$}(k+1) &= \text{LAT}(k),~
    \text{Raf}(k+1) = \text{Ras}(k)\\
    \text{Ras}(k+1) &= \text{Grb2Sos}(k) \vee \text{RasGRP1}(k)\\
    \text{RasGRP1}(k+1) &= \text{DAG}(k) \wedge \text{PKCth}(k)\\
    \text{Rlk}(k+1) &= \text{Lck}(k),~
    \text{Rsk}(k+1) = \text{ERK}(k)\\
    \text{SEK}(k+1) &= \text{PKCth}(k),~
    \text{SLP76}(k+1) = \text{Gads}(k)\\
    \text{TCR$^+$}(k+1) &= \neg \text{cCbl}(k) \wedge \text{TCRlig}(k)\\
    \text{TCR$^\dagger$}(k+1) &= \text{Fyn}(k) \vee [\text{Lck}(k) \wedge \text{TCR$^+$}(k)]\\
    \text{ZAP70}(k+1) &= \neg \text{cCbl}(k) \wedge \text{Lck}(k) \wedge \text{TCR$^\dagger$}(k)
    \end{align*}
    
\noindent where $(\cdot)^+$ represents the binding status, $(\cdot)^*$ denotes being activated, $(\cdot)^\dagger$ indicates phosphates,
CD45, CD4 and TCRlig are the three control inputs and there are four biomolecular species to be considered
as the outputs of the networks, which are AP1, CRE, NFAT and NFkB. 
In \cite{9969868}\footnote{The T-cell model in \cite{9969868} adopted from \cite{klamt2006methodology} 
includes a few typos.}, this T-cell Boolean network was used for structural controllability analysis to design a feedback controller. 

\emph{Case 1 ({\rm CD45=1, CD4=1, TCRlig=1}):} AP1 and NFAT converges to 0, NFkB converges to 1 and
CRE switches between the following update rules:
If $n$ is even and greater than or equal to 4,
\begin{align}
    \text{CRE}(k+n) &=
        \text{Lck}(k) \wedge \text{TCR$^\dagger$}(k) 
        \wedge \neg \text{cCbl}(k) \notag\\
        \wedge &[\text{PAGCsk}(k) \vee \text{ZAP70}(k) 
        \vee \neg \text{Fyn}(k)]\notag\\ 
        \wedge &[\text{PAGCsk}(k) \vee \text{ZAP70}(k) 
        \vee \neg \text{TCR$^+$}(k)]
\end{align}
and if $n$ is odd and greater than or equal to 5,
\begin{align}
    \text{CRE}(k+n) &=
        \text{cCbl}(k) \wedge \text{PAGCsk}(k) \wedge 
        \neg \text{ZAP70}(k) \notag\\
        &\wedge [\text{Fyn}(k) \vee \neg \text{TCR$^+$}(k)]
\end{align}
Based on this finding, an additional feedback control input would
be designed to derive CRE towards desired states.
A total of 21 states update logic out of the 37 states converge to a period of 2 switching logic.
The rest 16 states converge to steady states of either 0 or 1.
The state space shrinks from $2^{37}(137\text{ billion})$ to $2^{21} (2\text{ million})$, 
which is only 0.0015\% of the original size of the state space. 

\emph{Case 2 (at least one of {\rm CD45, CD4} or {\rm TCRlig} equal to 0):}
AP1, CRE and NFAT converge to 0 and NFkB converges to 1. 
There is no possibility of introducing further control structures to change
the outputs.

These analyses lead the problem space from originally computationally infeasible to feasible ranges.
In addition, they clearly show what states can or cannot be controlled and what the update-rule structures of states
to be controlled are.

From a biological perspective, when foreign antigens are presented to T cell receptors (\text{TCR}) and co-receptors (\text{CD4}), along with the involvement of receptor-type protein tyrosine phosphatase (\text{CD45}) \cite{klamt2006methodology}, this triggers an immune response signaling cascade that includes the Ras-Raf-Mek-Erk pathway.
This, in turn, leads to the transcriptional activation of numerous immune-related genes, such as IL-2, IL-6, IL-10, TNF-$\alpha$, and others \cite{wen2010role}.
The promoters of these genes commonly contain a DNA target sequence known as the cAMP-responsive element (\text{CRE}), to which transcription factors belonging to the CREB family (\text{CREB}) can specifically bind and initiate transcription.

To counterbalance the risk of an overactive immune response that could potentially result in autoimmunity, a negative feedback mechanism is in place.
This mechanism is primarily mediated by the E3 ubiquitin ligase (\text{cCbl}) \cite{duan2004cbl}, which facilitates the degradation of key signaling proteins, including the activated protein tyrosine kinase (\text{ZAP70}), thereby ensuring immune response homeostasis \cite{klamt2006methodology}.
Within the T-cell receptor network, with all control inputs set to 1, a negative feedback loop between \text{cCbl} and \text{ZAP70} continuously operates, resulting in oscillations between the states (\text{cCbl}=0, \text{ZAP70}=1) and (\text{cCbl}=1, \text{ZAP70}=0).
In the former states, according to the Boolean logic, \text{TCR$^+$}, \text{Fyn}, \text{PAGCsk}, \text{Lck}, and \text{TCR$^\dagger$}, which also form feedback regulations to \text{ZAP70}, attain values of 1, 1, 0, 1, and 1, respectively, while in the latter states, 0, 0, 1, 0, and 0, respectively.
In the former scenario, \text{ZAP70} gets activated by the feedback regulations and gives a positive signal to the downstream Ras-Raf-Mek-Erk pathway, leading to the activation of \text{CRE}.
Rather, in the latter scenario, \text{ZAP70} gets inactivated, so that \text{CRE} also switches back to the inactive state.
As a result, the long-term dynamics of \text{CRE} exhibit a switching behavior between 0 and 1, signifying that the T cell signalling cascade effectively orchestrates immune responses without straying into an excessive territory, thereby maintaining immune homeostasis.

This regulatory process corresponds to the results of \emph{Case 1 ({\rm CD45=1, CD4=1, TCRlig=1})}, where the recursive self-composite logic of \text{CRE} switches between (22) and (23).
In (22) and (23), Ras-Raf-Mek-Erk pathway between \text{CRE} and upstream feedback loops is omitted, while only the feedback regulations are represented as a switching for the presence or absence of the negation operator in front of \text{cCbl}(k), \text{ZAP70}(k), and \text{Fyn}(k).
In the results of \emph{Case 2}, on the other hand, when any one of the input controls equals to 0, multiple feedback loops will no longer be able to operate, leading to \text{ZAP70} keeping inactivated so that \text{CRE} converges to 0.
This signifies that the immune response signalling cascade is deactivated, rendering the homeostatic effect of negative feedback unnecessary.

In summary, we showed that the recursive self-composite logics of \text{CRE}, (22) and (23), serve as an intuitive representation of the long-term oscillatory dynamics of the output of the T cell Boolean network model.
\end{example}
\bigskip

\begin{example}[Cancer signaling network]
    In the study by Fumi\~a et al. \cite{fumia2013boolean}, a Boolean network model was developed to represent human tumorigenesis. This model encompasses 90 states of proteins within cancer signaling pathways, alongside 6 control inputs – Mutagen, GFs (growth factors), Nutrients, TNF$\alpha$, Hypoxia and Gli. The Boolean functions in the network were initially formulated using algebraic operators ($+$, $-$) and $\text{sgn}(\cdot)$ – a thresholding function in which $\text{sgn}(x)=0$ for $x\le 0$ and $1$ for $x>0$. As this algebraic representation is less straightforward for recursive self-composition, we opted to transform the Boolean functions into an equivalent form utilizing the Boolean operations, while preserving their original truth tables, the Boolean model is available
    to download as indicated in Supplementary Material. 
    
    The network outputs include two bio-molecules, Glut1 and Lactic acid, as well as two virtual nodes that represent cellular phenotypes – Apoptosis and DNA repair. In addition, the other cellular phenotype, Proliferation, was defined by examining the long-term dynamics of the network, particularly the activation sequence of cyclin nodes. In \cite{chu2015precritical}, attractor-transition analysis of this network was carried out to investigate the critical transition of tumorigenesis along with the accumulation of driver mutations.

    The update logic for Cyclins D has the longest string. Its update rule
    has 23k characters. 
    And, the first three self-compositions make the length increase exponentially, 
    i.e., 149k, 34M and 293M. Applying the self-composition
    procedure to the original network produces long string chains of
    the updated rules. The lengths of the strings are too long to be handled
    by most digital computers. As analyzing Boolean networks are known
    to be an NP-hard problem \cite{akutsu2007control}, 
    finding some networks producing large string beyond the current computer calculation speed and 
    memory capacity is not surprising.

    To restrict the network to the condition
    with normoxic (normal oxygen level), sufficient growth factors and nutrient-rich
    five input states are set as follows:
    Mutagen = 0, GFs = 1, Nutrients = 1, TNF$\alpha$ = 0, Hypoxia = 0 and Gli = 0. 
    In \cite{fumia2013boolean}, it is shown that the cell enters into the proliferation cycle with a period of 7.

    p53, Cyclins A and Cyclins D are identified by trial and error
    as the additional control inputs to change the phenotype. 
    Pinning the three states to 1 and applying the self-compositions make the logic length become too large again before it converges. 
    
    By applying the leaping every 4 iterations, i.e., starting from ${\bf x}(k+1)$,
    we obtain ${\bf x}(k+2)$, ${\bf x}(k+3)$, ${\bf x}(k+4)$ and ${\bf x}(k+5)$. Then, leaping
    from ${\bf x}(k+5)$, we obtain ${\bf x}(k+10)$, ${\bf x}(k+15)$ and ${\bf x}(k+20)$, and it converges at ${\bf x}(k+15)$. 
    And, applying the filling procedure from ${\bf x}(k+16)$ to ${\bf x}(k+19)$ we confirm the logic converges to a single logic.
    
    The maximum string length of the updated rules for each iteration is as follows: 
    17k, 29k, 60k, 65k, 64k, 10 and 10, where the corresponding longest string states are TSC1/2, Cytoc/APAF1, Cytoc/APAF1, GSH, eEF2k, E2F and E2F, respectively. 
    It converges to a steady state instead of the cycle and the phenotype changes from proliferation to apoptosis. More interestingly, all states converge to 0 or 1 regardless of the initial conditions except $\text{E2F}(k+15)$, a transcription factor for cell‐cycle regulation genes, which converges to
    \begin{align*}
        \text{E2F}(k+15) = \text{E2F}(k) \wedge [\neg \text{Rb}(k) ]
    \end{align*}
    where Rb is retinoblastoma protein. 
    
    The example demonstrates the power of the proposed method in finding a hidden simple logic behind the complex networks. Such finding can help to identify further important drug target candidates to be developed for efficient cancer therapeutic strategies. 
\end{example}

\section{Conclusions \& Future Works}
We present the recursive self-composite approach to reveal the hidden characteristics of Boolean networks. Most interestingly, we found the cyclic nature of synchronous Boolean update rules. This is the first time that the converging nature of the Boolean logic dynamics is unveiled explicitly. We also found several interesting properties of the converged logic: the existence of kernel logic and its relationship with the length of periodic cycles in the state space. There might be many other interesting hidden structures of the Boolean network, which we aim to reveal in future studies using the recursive self-composite approach. Finding fundamental relationships between the repeating logic and biological phenomena to control the behaviour of Boolean networks and extending the approaches to asynchronous Boolean networks are of immediate interest.

\section*{Supplementary Material}
 The Boolean network model of the cancer signalling network in Python is available to download at the following link: \href{https://github.com/myjr52/Fumia_cancer_network_boolean_model}{\url{https://github.com/myjr52/Fumia\_cancer\_network\_boolean\_model}}

\section*{Acknowledgment}
The authors would like to thank the support of the Cheney fellowship.
The research was initiated in June 2023 by the visit of KHC to the University of Leeds as one of the Cheney fellowship activities, hosted by JK. 

\ifCLASSOPTIONcaptionsoff
  \newpage
\fi

\bibliographystyle{IEEEtran}
\bibliography{IEEEabrv,recompBN}
\vfill

\end{document}